\documentclass{conm-p-l}

\newtheorem{theorem}{Theorem}[section]
\newtheorem{lemma}[theorem]{Lemma}
\theoremstyle{definition}
\newtheorem{definition}[theorem]{Definition}
\newtheorem{example}[theorem]{Example}
\newtheorem{proposition}[theorem]{Proposition}
\newtheorem{corollary}[theorem]{Corollary}

\theoremstyle{remark}
\newtheorem{remark}[theorem]{Remark}
\numberwithin{equation}{section}
\usepackage{latexsym}
\usepackage{amssymb}
\usepackage{epsfig}
\usepackage{graphics}

\newcommand{\bea}{\begin{eqnarray}}

\newcommand{\eea}{\end{eqnarray}}
\newcommand{\nn}{\nonumber \\}
\newcommand{\be}{\begin{equation}}

\newcommand{\ee}{\end{equation}}

\begin{document}
\title[Weak modules and logarithmic intertwining operators]{Weak
modules and logarithmic intertwining operators for vertex operator
algebras }
\author{Antun Milas}
\curraddr{Department of Mathematics, Rutgers University,
110 Frelinghuysen Road, Piscataway, NJ 08854-8019}
\email{amilas@math.rutgers.edu}
\thanks{The author is partially supported by NSF grants.}
\date{}
\bibliographystyle{amsalpha}

\begin{abstract}
We consider a class of weak modules for vertex operator algebras
that we call logarithmic modules. We also construct nontrivial examples
of intertwining operators between certain logarithmic modules for the Virasoro
vertex operator algebra. 
At the end we speculate about some possible
logarithmic intertwiners at the level $c=0$.
\end{abstract}

\maketitle

\section*{Introduction}

This work is an attempt to explain an algebraic reformulation of ``logarithmic
conformal field theory'' from the vertex operator algebra point of view.

The theory of vertex operator algebras, introduced in works of Borcherds
(cf. \cite{Bo}), Frenkel, Huang, Lepowsky and Meurman (\cite{FHL}, \cite{FLM})
, is the mathematical counterpart of 
conformal field theory, introduced in \cite{BPZ}.
One usually studies {\em rational} vertex operator algebras,
which satisfy a certain semisimplicity condition for modules.
If we want to go beyond rational vertex operator algebras we
encounter several difficulties. 
First, we have to study indecomposable, reducible  modules, for
which there is no classification theory. 
Another problem is that the notion of intertwining operator, 
as defined in \cite{FHL}, is not the most general.

We can illustrate this with the following example.
Let $L(c,0)$ (cf. \cite{FZ}) be a non--rational vertex operator algebra
associated to a lowest weight representation of the Virasoro
algebra and 
\begin{equation} \label{triple}
{\mathcal Y}_1, {\mathcal Y}_2
\end{equation}
a pair of intertwining operators for certain triples of $L(c,0)$--modules. 
It is well known 
that one can study matrix coefficients
\begin{equation} \label{matrix}
\langle w'_3, {\mathcal Y}_1(w_1,x)w_2 \rangle
\end{equation}
and
\begin{equation} \label{matrix1}
\langle {v'_4}, {\mathcal Y}_1(v_1,x_1){\mathcal Y}_2(v_2,x_2)  v_3 \rangle,
\end{equation}
where $w_i$, $v_j$, $i=1,2$, $j=1,2,3$, $v'_4$ and $w'_3$ are some vectors, 
by solving appropriate differential equations (see \cite{BPZ} or
\cite{H1}).
In several examples we find, after we switch to complex variables,
that in certain domains some solutions have logarithmic singularities. 
This fact was exploited in \cite{Gu1}, where
such conformal field theories are called ``logarithmic''.
The appearance of logarithms contradicts the following property
of intertwining operators:
\begin{equation} \label{prva}
{\mathcal Y}(w_1,x) \in {\rm Hom}(W_2,W_3)\{ x \},
\end{equation}
where $w_1 \in W_1$ and ${\rm Hom}(W_2,W_3) \{x\}$ denotes the
space of formal sums of the form $\sum_{n \in \mathbb{C}}v_n x^n$, 
$v_n \in {\rm Hom}(W_2,W_3)$. 
The aim of this paper is to develop an algebraic theory that 
incorporates logarithmic solutions as well.
We extend the property (\ref{prva}) (cf. Definition \ref{maindef}) so that we allow
\begin{equation} \label{druga}
{\mathcal Y}(u,x) \in {\rm Hom}(W_2,W_3)\{ x \}[{\rm log}(x)].
\end{equation}
Here ${\rm log}(x)$ is just another formal
variable, and {\em not}
a formal power series. Note that we allow arbitrary nonnegative
integral
powers of ${\rm log}(x)$ and arbitrary complex powers of $x$.
Also, we assume
the property $\frac{d}{dx}{\rm log}(x)=\frac{1}{x}$.
Such modification dictates, if we assume the usual condition that
${\mathcal Y}(L(-1)w,x)=\frac{d}{dx}{\mathcal Y}(w,x)$,
that the space $W_3$ is no 
longer an $L(c,0)$--module
but rather a weak module that is not $L(0)$--diagonalizable.
For this purpose we define a ``logarithmic module'' to be a weak module 
that admits a decomposition into generalized $L(0)$--eigenspaces.
The analogy with logarithms is clear 
when we consider ordinary differential equations 
(cf. Proposition \ref{structureint}).
If one carefully analyzes our definition
it is clear that logarithmic intertwiners can be defined solely
without logarithms, but then the formulas are not very transparent.

The next problem is to construct a non--trivial (not covered 
by the definition in \cite{FHL}) example of a logarithmic intertwining operator.
In this paper we study only logarithmic intertwining operators
associated to Virasoro vertex operator algebras.
By carefully examining Frenkel-Zhu-Li's formula (cf. \cite{FZ}, \cite{Li1}) 
we found sufficient
conditions (Theorem \ref{main}) for the existence 
of a non--trivial logarithmic intertwining operator each of the types
\begin{equation} \label{1log}
\binom{W}{L(c,h_1) \ M(c,h_2)} \  {\rm and} \  \binom{W}{L(c,h_1) \ L(c,h_2)},
\end{equation}
where $W$ is some 
where $W$ is some {\em logarithmic} module.
Then we apply Theorem \ref{main} for
some special theories: $c=-2$ and $c=1$ (Corollary \ref{c=1} and
Corollary \ref{c=-2}).
Note that we limit ourselves to the case when
$W_1$ and $W_2$ are ordinary modules.
It seems to be harder to construct 
a logarithmic intertwining operator of the type
$\binom{W_3}{W_1 \ W_2}$ 
where $W_i$, $i=1,2$ are non--trivial logarithmic modules.

It is clear that non--trivial logarithmic modules do not
come up in the theory of rational 
vertex operator algebras; therefore they do not come up 
in the case of the vertex operator algebras
$L(c_{p,q},0)$, where $c_{p,q}=1-\frac{6(p-q)^2}{pq}$, $p,q \geq 2$
and $(p,q)=1$. 
In particular, if $c=c_{2,3}=0$ then
$L(0,0)$ is the trivial vertex operator algebra (hence
the representation theory is trivial). 
But if we consider the vertex operator algebra 
$M_0=M(0,0)/ \langle L(-1)v_{0} \rangle $
(here $v_0$ is the lowest weight vector)
the situation is substantially different.
We can consider modules for
$M_0$ and study the corresponding intertwining
operators. What is interesting is that $M_0$ is not a simple vertex
operator algebra and this makes the whole subject very interesting.
We should mention that the $c=0$ case has been studied in the 
connection with {\em percolation} in mathematics and physics
(for a good review of the subject see \cite{LPS} and the references therein) 

Our motivation to study the $c=0$ case stems from
\cite{GL} where Gurarie and Ludwig studied logarithmic
OPE (operator product expansions) between certain special $c=0$
theories.
Their result, slightly modified in our language,
predicts logarithmic intertwiners
of the type
\begin{equation} \label{unknown}
\binom{W_b}{L(0,\frac{m^2-1}{24})  \ L(0,\frac{m^2-1}{24})},
\end{equation}
for a certain logarithmic module $W_b$. 
The constant $b$, as shown in \cite{GL},
seems to be the same for various $m$. 
It is unclear in \cite{GL} what the vector space
$W_b$ is and how to construct intertwining operators of the 
type (\ref{unknown}). 

We speculate what properties the space $W_b$ should have, and explain
why the constant $b$ appears instead of the central charge. 
Consider the Lie
algebra $\hat{W}_{{\rm log}}$ \footnote{A related algebraic structure  
was considered in \cite{GL}.}
, generated by $t^{(i)}(n)$, $i,n \in \mathbb{Z}$,
$i \in \mathbb{Z}$ and $b$, with  commutation relations
$$[t^{(i)}(m),t^{(j)}(n)]=(m-n)t^{(i+j)}(m+n)+(j-i)t^{(i+j-1)}(m+n)+
{\rm c}(t^{(i)}(m),t^{(j)}(n))b,$$
where ${\rm c}(t^{(i)}(m),t^{(j)}(n))$ is a certain 2--cocycle introduced
in the paper and $b$ is the central element.
Note that thera are two (distinguished) 
Virasoro subalgebras inside $\hat{W}_{{\rm log}}$. 
{\em Horizontal}, take $i=j=0$ and {\em vertical}, take $m=n=0$. 
It will turn out that both the horizontal and vertical Virasoro subalgebras 
have central charge $0$ no matter what is $b$. 
We show that there is an evidence that logarithmic intertwining
operators exist but it is much more difficult to construct 
them and explicitly define what is the space $W_b$.
In a sequel we will study this problem in more details. 
Also, we believe that studying the representation theory of $\hat{W}_{{\rm log}}$ is an interesting subject. \\
{\em n.b.} A lot of progress has been done by physicists 
in understanding LCFT.  Instead of listing numerous references we refer an interested reader to the excellent review articles \cite{G}, \cite{F1} and \cite{RT} (and references therein). \\
{\bf Acknowledgment:} I benefited
from discussions with V. Gurarie, Y.-Z. Huang and J. Lepowsky.

\section{Weak modules and logarithmic modules}

\subsection{Notation}

We denote by $\mathbb{N}$ the
set of positive integers, by $\mathbb{Z}$ the ring of integers
and by $\mathbb{C}$ the field of complex numbers.
In what follows all variables are formal. For a vector space $W$ 
we denote by $W \{x \}$
the vector space of formal sums $\sum_{r \in
\mathbb{C}} a_r x^r$, where $a_r \in W$, 
by $W[[x.x^{-1}]]$ the vector space 
of formal Laurent series and by $W((x))$ the
vector space of truncated Laurent series. In particular, $\mathbb{C}((x))$
is the ring of truncated Laurent series. 

\subsection{Logarithmic modules}

Let $V=\bigoplus_{n \in \mathbb{Z}} V_n$ 
be a vertex operator algebra (for the definition see \cite{FHL} or \cite{FLM}).
A {\em weak} $V$--module $W$ by definition (cf. \cite{DLM}) satisfies all
the axioms for a $V$--module (see \cite{FHL}) 
except for those involving the grading and the action of
$L(0)$. In particular
we do not assume any grading on $W$. Let $\mathbb{I} \subset \mathbb{C}$.
We say that a weak $V$--module $W$ is $\mathbb{I}$--gradable (cf. \cite{Zh1})
if there exists a decomposition 
$$W=\bigoplus_{i \in \mathbb{I}} W(i),$$
such that
$$v_n W(i) \subset W(i+{\rm deg}(v)-n-1),$$
for $v \in V({\rm deg}(v))$, where ${\rm deg}(v)$ is the degree
of $v$ (cf. \cite{FLM}) and $n \in \mathbb{Z}$. 
If every $\mathbb{N}$--gradable weak $V$--module is completely
reducible (in the category of weak $V$--modules) then we say that 
$V$ is rational \footnote{The reader should be aware that there are
several distinct 
notions of rationality in the literature.}. 
In this paper most of the statements 
are vacuous if $V$ is rational.

Let us consider a special class of gradable weak modules. 
\begin{definition} \label{logarithm}
We say that a weak $V$--module $W$ is {\em logarithmic}
if it admits a direct sum decomposition into 
generalized $L(0)$--eigenspaces, i.e. admits a Jordan form with
respect to the action of $L(0)$.
\end{definition}
Let $W$ be a logarithmic module.
Put 
\begin{equation} \label{1}
W_k:=\{ w \in W : (L(0)-k)^n w=0, {\rm for} \ {\rm some} \ n \in \mathbb{N} \}.
\end{equation}
\begin{proposition}
Let $W$ be as above and
$W=\bigoplus_{k \in \mathbb{N}+h} W_k$ for some $h \in \mathbb{C}$.
Then $W$ is an $\mathbb{N}$--gradable module.
\end{proposition}
\begin{proof}
Suppose that $v \in V(l)$ and $w \in W_k$. From the Jacobi identity
it follows that
$$(L(0)-(l+k-m-1))v_m w=v_m (L(0)-k)w,$$
for every $m$.
Therefore
$$(L(0)-(l+k-m-1))^n v_m w=0,$$
Thus $v_m w \in W_{l+k-m-1}$.
If we let
$$W(k):=W_{k+h},$$
then $W=\bigoplus_{k \geq 0} W(k)$ is $\mathbb{N}$--gradable.
\end{proof}

Let ${\mathcal C}_n$ be the category \footnote{A related category
$\mathcal O$ has
been considered in \cite{DLM}.} whose objects are 
logarithmic $V$--modules such that in (\ref{1}) 
the constant $n \in \mathbb{N}$ is fixed, i.e., with respect to
the action of $L(0)$, $W$ has a matrix realization with Jordan blocks
of size at most $n$. 
Then we have  
$${\mathcal C}_1 < {\mathcal C}_2 < {\mathcal C}_3
< \ldots .$$
In this paper we shall work mostly with the category ${\mathcal C}_2$.
We say that a logarithmic module $W$ is nontrivial
if it  is not an object in ${\mathcal C}_1$.
Note that ${\mathcal C}_1$ is the category of
all weak $V$--modules that are
$L(0)$--diagonalizable and $\mathbb{N}$--gradable.
Often, we are not interested in studying logarithmic modules  
because every irreducible logarithmic 
$V$--module is $L(0)$--diagonalizable,
i.e., it is an object in ${\mathcal C}_1$.

\subsection{Intertwining operators for logarithmic modules}

Here we extend
the definition of intertwining operators for a triple of modules as stated in \cite{FHL}.

Let $\mathbb{C}[[x,x^{-1}]][y]$ be the space of formal Laurent series
with coefficients in $y$. Let 
$\frac{d}{dx} \in {\rm End}(\mathbb{C}[[x,x^{-1}]]).$
We extend $\frac{d}{dx}$ to $\mathbb{C}[[x,x^{-1}]][y]$ 
by letting $\frac{d}{dx}y=\frac{1}{x}$.

For obvious reasons we will write ${\rm log}(x)$ instead of $y$, where
${\rm log }(x)$ is purely formal and {\rm not} a power series.
If $A(x) \in x\mathbb{C}[[x]],$
we set
\begin{equation} \label{log1}
{\rm log}(1+A(x))=\sum_{n >0} (-1)^{n-1} \frac{A(x)^n}{n} \in \mathbb{C}[[x]].
\end{equation}
\begin{definition} \label{maindef}
Let $W_1$, $W_2$ and $W_3$ be three logarithmic $V$--modules.
A logarithmic intertwining operator 
is a linear mapping
\bea
&& {\mathcal Y}( \ ,x) : W_1 \rightarrow {\rm Hom}(W_2,W_3) \{x\}[{\rm log}(x)] \\
&& {\mathcal Y}(w,x)=\sum_{i \geq 0} \sum_{ \alpha \in \mathbb{C}}
w^{(i)}_{\alpha}x^{-\alpha-1} {\rm log}^i(x),
\eea
such that the following properties hold:

\begin{enumerate}

\item The {\it truncation} property: For any $w_i \in W_i$, $i=1,2$,  
and $j \geq 0$
$$(w_1)^{(j)}_\alpha w_2=0,$$
for large enough ${\mathcal Re}(\alpha)$.

\item The {\it $L(-1)$-derivative property}: For any $w_1\in W_1$, 
$$\mathcal{Y}(L(-1) w_1, x)=\frac{\partial}{\partial x}
\mathcal{Y}(w_1, x).$$

\item The {\it Jacobi identity}: In $W_{3} \{x_{0},x_{1},x_{2}\} [{\rm log}(x_2)]$, we have 
\bea \label{jac}
\lefteqn{x_0^{-1} \delta \left ( \frac {x_1-x_2}{x_0} \right ) 
Y(u, x_1) \mathcal{Y}(w_1, x_2)w_2} \\
&& -x_0^{-1} \delta \left ( \frac
{x_2-x_1}{-x_0} \right ) \mathcal{Y}(w_1, x_2) 
Y(u, x_1)w_2=
\nn
&&= x_2^{-1} \delta \left ( \frac {x_1-x_0}{x_2} \right
)\mathcal{Y}(Y(u, x_0)w_1, x_2)w_2, \nonumber
\eea
for $u\in V$, $w_1 \in W_1$ and $w_2 \in W_2$.
\end{enumerate}
\end{definition}

In this paper
we denote the space of all logarithmic intertwining operators of the type
$\binom{W_3}{W_1 \  W_2}$
by $I \binom{W_3}{W_1 \ W_2}$ \footnote{Notice that every intertwining
operator (as defined in \cite{FHL}) is also a logarithmic 
intertwining operator. 
We say that a logarithmic operator is {\em non--trivial} if it is not covered
by the definition in \cite{FHL}.}.

\begin{remark} 
Let 
$${\mathcal Y}^{(i)}(w,x)=\sum_{n \in \mathbb{C}} w^{(i)}_n x^{-n-1}.$$
It follows from (\ref{jac}), by extracting
the coefficients of ${\rm log}^i(x)$, that  
Jacobi identity holds for each ${\mathcal Y}^{(i)}$ separately. Hence 
our definition can be stated without logarithms, but then the
$L(-1)$--derivative property does not hold for each ${\mathcal Y}^{(i)}$.
\end{remark}
Here are some consequences of Definition \ref{maindef}.
\begin{proposition}
Let ${\mathcal Y} \in I  \binom{W_3}{W_1 \ W_2}$ and $w \in W_1$.
Then
\begin{itemize} \item[$(a)$]
$$e^{yL(-1)}{\mathcal Y}(w,x)e^{-yL(-1)}={\mathcal Y}(w,x+y),$$
where
$${\rm log}(x+y):=e^{y \frac{d}{dx}}{\rm log}(x).$$
\item[$(b)$]
$$e^{yL(0)}{\mathcal Y}(w,x)e^{-yL(0)}={\mathcal Y}(e^{y L(0)} w,e^yx),$$
where ${\rm log}(e^yx)=y+{\rm log}(x)$. 
\item[$(c)$]
Define 
$${\mathcal Y}^*(w_2,x)w_1=e^{xL(-1)}{\mathcal Y}(w_1,e^{\pi i} x)w_2.$$
Then ${\mathcal Y}^* \in I \binom{W_3}{W_2 \ W_1}$
where ${\rm log}(e^{\pi i}x):={\rm log}(x)+ \pi i$.
\end{itemize}
\end{proposition}
\begin{proof}
(a) follows from the $L(-1)$--derivative property and
$$[L(-1),{\mathcal Y}(w_1,x)]=\frac{d}{dx}{\mathcal Y}(w_1,x).$$
For (b) we use the same proof as in \cite{FHL} Lemma 5.2.3, i.e.
formulas $$[L(0),{\mathcal Y}(w_1,x)]={\mathcal Y}(L(0)w_1,x)+x{\mathcal Y}(L(-1)w_1,x)$$
and
$$e^{y x \frac{d}{dx}} {\rm log}(x)=y+{\rm log}(x).$$ 
The statement (c) can be easily checked by using the same 
proof as in \cite{FHL}.
\end{proof}
\subsection{Some representations of the Virasoro algebra} \label{alstru}

\begin{definition}
We denote by ${\rm Vir}$ the Virasoro Lie algebra with the standard
triangular decomposition
${\rm Vir}={\rm Vir}_+ \oplus {\rm Vir}_0 \oplus {\rm Vir}_-.$ 
Let $M$ be a ${\rm Vir}$--module. We say that $M$ is a restricted
$Vir$--module if for every $v \in M$, $L(n)v=0$ for $n>>0$. Denote by
${\mathcal R}_c$ the category of all restricted ${\rm Vir}$--modules
of central charge $c$.
\end{definition}

Let us denote by $M(c,h)$ the Verma module for the Virasoro Lie
algebra of central charge $c$ and lowest weight $h$ 
and by $L(c,h)$ the corresponding irreducible quotient.
We write ${\mathcal O}_c$ for the category generated by all
lowest weight modules
of the central charge $c$.
Objects in ${\mathcal O}_c$ are restricted $L(0)$--diagonalizable 
${\rm Vir}$--modules such that for every $W$ in ${\mathcal O}_c$
\begin{equation} \label{spectrum}
{\rm Spec}L(0)|_W \in h_i+\mathbb{N},
\end{equation}
for some $h_i \in \mathbb{C}$, $i=1,...,k$.

Motivated by the definition of logarithmic modules 
we enlarge the category ${\mathcal O}_c$ in the following
way. 
Denote by ${\mathcal O}_c^n$ the category of all restricted $Vir$--modules
that admit decompositions into generalized $L(0)$--eigenspaces
with Jordan blocks of size at most $n$. 
Again all irreducible objects of ${\mathcal O}_c^n$ are
contained in ${\mathcal O}_c$.
For every $c,h \in \mathbb{C}$ and $n \in \mathbb{N}$ 
we define a {\em generalized} Verma module
$$M_{n}(c,h)={\rm Ind}_{U({\rm Vir}_{ \geq 0})}^{U({\rm Vir})}
{\mathcal V}(n),$$
where ${\mathcal V}(n)$ is a  $n$--dimensional $U({\rm Vir}_+) \oplus
\mathbb{C}L(0)$--module \footnote{More generally,  one can take an arbitrary
finite--dimensional $\mathbb{C}[L(0)]$--module.} 
 with a basis $v_i$, $i=1,...,n$ such that
\bea
&& L(n).v_i=0, n >0, \\
&& L(0).v_1=hv_1, \nn
&& L(0).v_i=hv_i+v_{i-1}, \ i=2,...,n,\nn
&& c|_{{\mathcal V}(n)}=c{\rm Id}. \nonumber
\eea
Clearly, $M_{n}(c,h)$ is an object in ${\mathcal O}_c^n$.
In the case $n=2$ we shall write $v:=v_1$ and $w:=v_2$.
Obviously for every $n \geq 2$ we have the following exact sequence:
\begin{equation} \label{dva}
0 \rightarrow M_{n-1}(c,h) \rightarrow M_n(c,h) \rightarrow M(c,h)
\rightarrow 0.
\end{equation} 
As before we will use the following notation
\be \label{type}
M_n(c,h)(m)=\{ v : (L(0)-m-h)^n v=0 \};
\end{equation}
$$M_n(c,h)=\bigoplus_{m \geq 0} M_n(c,h)(m).$$

Let $M$ be a $Vir$--module.
We say that $u \in M$ is a primitive vector (cf. \cite{K}) if there is
a submodule $U \subset M$ such that $U({\rm Vir}_+) u \in U$.
A primitive vector $u$ is a singular vector if $U=0$.
The most important property of Verma modules for the Virasoro
algebra is that every submodule is generated by its singular vectors. 
For the modules $M_n(c,h)$, $n \geq 2$ this is no longer true.
The structure of $M_n(c,h)$ is in general more complicated
then the structure of the ordinary Verma modules (cf. \cite{FF1}).

\begin{proposition} \label{primitive}
Every submodule of $M_2(c,h)$ is generated by its primitive
vectors. In general there are submodules
that are not generated by its singular vectors.
\end{proposition}
\begin{proof}
The first statement is clear.
Consider a Verma module $M_2(0,0)$
and a submodule $M'$ generated by $v$ and $L(-1)w$, i.e.
$M'$ is characterized by a non--split exact sequence
$$0 \rightarrow M(0,0) \rightarrow M' \rightarrow M(0,1) \rightarrow
0.$$
If $M'$ is generated by its singular vectors then it has to be 
generated by $v$ and some singular vector belonging
to $M_2(0,0)(1)$. Hence it is of the form
$\alpha L(-1)v+\beta L(-1)w$, $\alpha,\beta \in \mathbb{C}$. It
is easy to check that any such singular vector is a multiple of
$L(-1)v$ (because of $L(1)L(-1)w=2v$). This would imply $L(-1)w \notin M'$.
Therefore $M'$ is not
generated by its singular vectors.
\end{proof}

Motivated by the notion of maximal module for
ordinary Verma modules we define  $M'_2(c,h) \subset M_2(c,h)$ 
as the union of all submodules $M'$ of $M_2(c,h)$ 
such that $M' \cap M_2(c,h)(0)=0$. 
Then $\tilde{W}_2(c,h):=M_2(c,h)/M'_2(c,h)$ 
does not contain singular vectors of the weight strictly bigger then $h$.
For computational purposes we will need a slighlty ``bigger''
representations (and their contragradient versions)
then $\tilde{W}_2(c,h)$. The following construction
will turn out to be useful.
Let $M_2(c,h)/M''$ be a $Vir$--module 
where $M''$ is the maximal $L(0)$--diagonalizable submodule 
such that $M'' \cup M_2(c,h)(0)=0$. Clearly $M^0_2(c,h) \subset M'_2(c,h)$.
We let $W_2(c,h)=(M_2(c,h)/M^0_2(c,h))'$ where $( \ )'$ stands for
the contragradient functor.
\begin{definition} \label{almostdef}
We say that a module $W_2 \in {\rm Ob}{\mathcal O}_c^2$ is {\em almost
irreducible} if it is isomorphic to 
$\tilde{W}_2(c,h)$ for some $c,h \in \mathbb{C}$.
\end{definition}
We have an exact sequence of $Vir$--modules:
\begin{equation} \label{almost}
0 \rightarrow L(c,h) \rightarrow W_2 \rightarrow \tilde{M}(c,h) \rightarrow 0,
\end{equation}
where $\tilde{M}(c,h)$ is some quotient of $M(c,h)$.
Especially interesting are the modules $\tilde{W}_2(c,h)$ 
such that there is a non--split exact sequence:
$$0 \rightarrow L(c,h) \rightarrow W_2(c,h) \rightarrow L(c,h) \rightarrow 0,$$
i.e., a non--trivial vector in ${\rm Ext}^1(L(c,h),L(c,h))$.
On the other hand
$$0 \rightarrow L(c,h) \rightarrow M_2(c,h)/M_2^0(c,h) \rightarrow M(c,h) \rightarrow 0$$
and
$$0 \rightarrow M(c,h)' \rightarrow W_2(c,h) \rightarrow L(c,h) \rightarrow 0
.$$
Because $L(c,h) \hookrightarrow M(c,h)'$, an important feature of
$W_2(c,h)$ is that irreducible module $L(c,h)$ both embedds and
is a quotient in
 $W_2(c,h)$.
Let us denote by $S_{M_n(c,h)}(m)$ the dimension of the singular
subspace for \\
$M_n(c,h)(m)$.
\begin{lemma} \label{primitive1}
For every $n \in \mathbb{N}$, ${\rm dim} \  S_{M_2(c,h)}(n) \leq 2$.
\end{lemma}
\begin{proof}
Let $v_{sing} \in M_2(c,h)$. From (\ref{dva}) it follows that
$v_{sing}=w_1.v+w_2.w$, $w_1,w_2 \in U(Vir_-)$, where $w_2.v$ is a singular
vector in $M(c,h)$.
If ${\rm dim}\ S(n) >2$ then this would imply that
there are three linearly independent singular vectors
of the form $w_1'.v+w_2.w$, $w_1.v$ and $w_1.v+w_2.w$ inside
$M_2(c,h)(n)$. Therefore
$(w_1'-w_1).v$ should be singular too. On the other hand,  
the singular space for $M(c,h)$ is at most one--dimensional 
for every $n$. Hence
$(w_1'-w_1).v$ is proportional to $w_2.v$. Therefore $w_1'.v+w_2.w$
is a linear combination of $w_2.v$ and $w_1.v+w_2.w$.
\end{proof}
More generally, by using the induction,
one can show that $S_{M_n(c,h)}(m) \leq n$ for every $m,n \in \mathbb{N}$.
%

\subsection{Generalized density modules}

Another important class of modules for the Virasoro algebra
are the so called density modules ${\mathcal F}_{\lambda,\mu}$ (cf. \cite{FF1},\cite{FF2}).
As a vector space ${\mathcal F}_{\lambda,\mu}$
can be realized as $x^{\mu} \mathbb{C}[x,x^{-1}] dx^{\lambda}$
(cf. \cite{FF2})
such that the generators $L_n$ act as certain vector fields.
This can be generalized by considering the space:
$$\bigoplus_{0 \leq i \leq n} x^{\mu}\mathbb{C}[x,x^{-1}]{\rm
log}^i(x)dx^{\lambda}.$$
Let
$${\mathcal F}_{\lambda,\mu,n,\beta}=\bigoplus_{0 \leq i \leq n, r \in 
\mathbb{Z}} \mathbb{C} u^{(i)}_{r},$$
such that 
$$L_m.u_{r}^{(i)}=(\mu+r+\lambda(m+1))u_{r-m}^{(i)}+ \beta iu_{r-m}^{(i-1)}.$$
Then ${\mathcal F}_{\lambda,\mu,n,\beta}$ is a $Vir$--module. Again
this module is reducible for $n \geq 1$. \\
{\em n.b} These modules (at least when $n=1$) are studied in \cite{CKW}
in connection with the classification of all extensions of
the {\em conformal modules} for (super)conformal algebras.

\subsection{Some properties of logarithmic intertwiners}

Suppose that  $V$ is a vertex operator algebra, $W_i$, $i=1,2$ are 
objects in ${\mathcal C}_1$ and $W_3$ is an object in ${\mathcal C}_k$.
Let $\mathcal Y$ be a logarithmic intertwining operator of the type
$ \binom{W_3}{W_1 \ W_2}$.
Assume that
${\rm Spec}L(0)|_{W_i} \in h_i + \mathbb{N}$.
Then for every  $w_1 \in W_1$, $w_2 \in W_2$ and $w'_3 \in W'_3$
such that $w_i \in W_i(m_i)$ for $i=1,2$ and
$w'_3 \in W'_3(n)$,
we have
\bea
&& \langle (L(0)-h_3-n)^k w'_3, {\mathcal Y}(w_1,x)w_2 \rangle= \nn
&& (x\frac{d}{dx}-h_3+h_1+h_2+m_1+m_2-n)^k \langle w'_3, {\mathcal Y}(w_1,x)w_2\rangle=0.
\eea
Therefore
\begin{equation} \label{onelog}
{\mathcal Y}(w_1,x)w_2 \in x^{h_3-h_1-h_2} W_3((x)) \oplus 
\cdots \oplus x^{h_3-h_1-h_2} W_3 ((x)) {\rm log}^{k-1}(x),
\end{equation}
for every $w_1 \in W_1$ and $w_2 \in W_2$.

Even a stronger statement holds.
\begin{proposition} \label{structureint}
Let $W_1 \in {\mathcal C}_{k_1}$, $W_2 \in {\mathcal C}_{k_2}$ and
$W_3 \in {\mathcal C}_{k_3}$ be logarithmic $V$--modules such that
${\rm Spec}(L(0))|_{W_i} \in h_i +\mathbb{N}$ for $i=1,2,3$
and 
$${\mathcal Y} \in I \ \binom{W_3}{W_1 \ W_2}.$$ Then
for every $w_1 \in W_1$, $w_2 \in W_2$
\bea \label{2int}
&& {\mathcal Y}(w_1,x)w_2 \in \nn
&& x^{h_3-h_1-h_2} W_3((x)) \oplus \ldots 
\oplus x^{h_3-h_1-h_2} {\rm log}^{k_1+k_2+k_3-3}(x)W_3((x)).
\eea
\end{proposition}
\begin{proof}
We may assume that $w_1$ and $w_2$ are homogeneous, i.e. they
are contained in some generalized eignespaces. Clearly
$(L(0)-h_1-m_1)^{k'_1}w_1=0$ and $(L(0)-h_2-m_2)^{k'_2}w_2=0$ for
some $h_1,h_2 \in \mathbb{C}$, $m_1,m_2 \in \mathbb{N}$ 
and $k'_1 \leq k_1$, $k'_2 \leq k_2$.
We will prove by induction on $k'_1+k'_2$ that
\bea \label{3int}
&& {\mathcal Y}(w_1,x)w_2 \in \nn
&& x^{h_3-h_1-h_2} W_3((x)) \oplus \ldots 
\oplus x^{h_3-h_1-h_2} {\rm log}^{k'_1+k'_2+k_3-3}(x)W_3((x)),
\eea
which implies the desired result. 
We already proved that if $k'_1=k'_2=1$ the statement holds for 
every $k_3$. Let $k'_1+k'_2>2$.
Suppose that the statement holds for every $k'_1$, $k'_2$ such that $k'_1+k'_2<k$. 
For every homogeneous $w'_3 \in W'_3$ there is $n$ such that 
$(L(0)-h_3-n)^{k_3}w'_3=0$. Then
\bea \label{genbeh}
&& \langle (L(0)-h_3-n)^{k_3}w'_3,{\mathcal Y}(w_1,x)w_2 \rangle= \nn
&& (x\frac{d}{dx}-h_3-n+L(0) \otimes 1 + 1 \otimes L(0))^{k_3} \langle
w'_3,{\mathcal Y}(w_1,x)w_2 \rangle=0,
\eea 
where 
$$1 \otimes L(0). {\mathcal Y}(w_1,x)w_2={\mathcal Y}(w_1,x)L(0)w_2$$
and
$$L(0) \otimes 1. {\mathcal Y}(w_1,x)w_2={\mathcal Y}(L(0)w_1,x)w_2.$$
The formula (\ref{genbeh}) can be written as 
\bea \label{gb1}
&& \biggl(x\frac{d}{dx}-h_3+h_1+h_2-n+m_1+m_2(L(0)-h_1-m_1 )\otimes 1 + \nn 
&& 1 \otimes (L(0)-h_2-m_2) \biggr)^{k_3} 
\langle w'_3,{\mathcal Y}(w_1,x)w_2 \rangle=0.
\eea
After expanding (\ref{gb1})---by using the binomial theorem---we obtain
\bea \label{gb2}
&& (x\frac{d}{dx}-h_3+h_1+h_2-n+m_1+m_2)^{k_3} \langle w'_3,{\mathcal
Y}(w_1,x)w_2 \rangle= \\
&& -\sum_{n \geq 1}^{k_3} \binom{k_3}{n}
(x\frac{d}{dx}-h_3+h_1+h_2-n+m_1+m_2)^{k_3-n} \cdot \nn 
&& ((L(0)-h_1 )\otimes 1 + 1 \otimes (L(0)-h_2))^n 
\langle w'_3,{\mathcal Y}(w_1,x)w_2 \rangle=0.\nonumber
\eea
If we apply the induction hypothesis, 
(\ref{gb2}) reduces to the following differential
equation
\bea \label{gb3}
&& (x\frac{d}{dx}-h_3+h_1+h_2+m_1+m_2-n)^{k_3} \langle w'_3,{\mathcal
Y}(w_1,x)w_2 \rangle = \nn
&&= x^{h_3-h_1-h_2+n-m_1-m_2}P({\rm log}(x)),
\eea
where ${\rm deg}(P) \leq k_3+k'_1+k'_2-4$.  Every solution
of (\ref{gb3}) inside $\mathbb{C}\{x \}[{\rm log}(x)]$
is of the form $$x^{h_3-h_1-h_2+n-m_1-m_2}Q({\rm log}(x)),$$
where $${\rm deg}(Q) \leq {\rm max}(k_3-1,k_3+k'_1+k'_2-4)+1
=k_3+k'_1+k'_2-3.$$
\end{proof}

\subsection{Vertex operator algebra $L(c,0)$, $c \neq c_{p,q}$}

Suppose that $c \neq
c_{p,q}:=1-\frac{6(p-q)^2}{pq}$ for every 
$p,q \geq 2$, $p,q \in \mathbb{N}$ and $(p,q)=1$. It is known that
the corresponding quotient $L(c,0)=M(c,0)/ \langle L(-1)v \rangle$ 
is irreducible.
$L(c,0)$ can be equipped with a structure of vertex operator algebra (cf. \cite{FZ}) and
all restricted ${\rm Vir}$--modules are weak $L(c,0)$--modules.
However, every irreducible $L(c,0)$--modules is of the form
$L(c,h)$, for some $h \in \mathbb{C}$. The category of all weak
$L(c,0)$--modules is too big for our purposes though.
Hence, we consider only the logarithmic modules in the
category ${\mathcal O}_2^c$. 
One of the main problems in vertex operator algebra theory is to
calculate fusion rules for a triple of modules. In the rational
vertex operator algebra setting it is enough to consider 
fusion rules for triples of
irreducible modules. Frenkel-Zhu's formula (\cite{FZ}, \cite{Li1})
is a valuable tool for such computation even in the non--rational case. 

Let us recall how to compute the dimension of the space
\begin{equation} \label{11m}
I \binom{L(c,h)}{L(c,h_1) \ L(c,h_2)}.
\end{equation}
for a triple of irreducible
$L(c,0)$--modules (cf. \cite{M1}, \cite{W}).

>From \cite{FZ} it follows that Zhu's algebra (\cite{Zh1},
\cite{FZ}) $A(L(c,0))$ is isomorphic to a polynomial algebra $\mathbb{C}[y]$.
Then for every  $L(c,0)$--module
$L(c,h)$ we can associate an $A(L(c,0))$--bimodule   
$A(L(c,h))$ (cf. \cite{FZ}, for details in our
setting see \cite{M1}, \cite{W}). It follows that
$$A(M(c,h)) \cong \mathbb{C}[x,y],$$
where 
$$y=[L(-2)-L(-1)], x=[L(-2)-2L(-1)+L(0)].$$
We will write
$$A(W_1,W_2):=A(W_1) \otimes_{A(V)} W_2 (0),$$
when there is no confusion about what is $V$.
By using results from \cite{FF1} (cf. \cite{M1})  
$$A(L(c,h_1),L(c,h_2)) \cong
\frac{\mathbb{C}[x]}{\langle p(x) \rangle},$$
for a certain polynomial $p(x)$. 
Then by using a result from \cite{FZ} the space (\ref{11m}) is $0$ 
if and only if  
$p(h) \neq 0$. If $c=c_{p,q}$, then $p(h)=0$ is 
also a sufficient condition for (\ref{11m}) to be one--dimensional. 
In general (cf. \cite{Li1}) if $p(h)=0$ we obtain 
a non--trivial intertwining operator of the type
$$\binom{L(c,h)}{L(c,h_1) \ M(c,h_2)}.$$ 
To obtain an intertwining operator of the
type $\binom{L(c,h)}{L(c,h_1) \ L(c,h_2)},$
it is enough to check that (cf. \cite{M1}) 
$${\rm Hom}_{A(L(c,0))}(A(L(c,h_2),L(c,h_1)),
L(c,h)(0)),$$
is non--trivial, i.e., that $q(h)=0$, where
$$A(L(c,h_2),L(c,h_1)) \cong
\frac{\mathbb{C}[x]}{\langle q(x) \rangle}.$$

\subsection{Main theorem}

In what follows we will give a sufficient condition
for the existence of a non--trivial logarithmic intertwining
operator of the type
\begin{equation} \label{typea}
\binom{M_2(c,h)'}{L(c,h_1) \ L(c,h_2)},
\end{equation}
Since the space $I \binom{M_2(c,h)'}{L(c,h_1) \ L(c,h_2)}$
is too big for our purposes 
($M_2(c,h)$ might have many submodules ) it is more
convenient to study the space
\begin{equation} \label{typeb}
I \binom{W_2(c,h)}{L(c,h_1) \ L(c,h_2)},
\end{equation}
where $W_2(c,h)$ is defined in the Chapter 1.

We shall try to follow a similar procedure as in the case of irreducible
modules (as explained in the previous section)
by carefully building logarithmic intertwining operators.
First, we obtain an upper bound for the
dimension of the space (\ref{typeb}).

Let $v,w \in M_2(c,h)(0)$ as before 
($L(0).v=hv$, and $L(0).w=hw+v$) and let
${\mathcal Y} \in I \binom{W_2(c,h)}{L(c,h_1) \ L(c,h_2)}.$
Define 
$$o^{(0)}_{\mathcal Y}(w_1)={\rm coeff}_{x^{h-h_1-h_2}}{\mathcal Y}(w_1,x)$$
and 
$$o^{(1)}_{\mathcal Y}(w_1)={\rm coeff}_{{\rm log}(x)x^{h-h_1-h_2}}{\mathcal
Y}(w_1,x).$$ Then $o^{(i)}_{\mathcal Y}$ defines
a linear map from $L(c,h_1)$ to ${\rm Hom}(L(c,h_2)(0),W_2(c,h)(0))$.
\begin{lemma} \label{l0}
For every $i$
the mapping $$\pi^i_{{\mathcal Y}} : A(L(c,h_1),L(c,h_2))
\rightarrow 
W_2(c,h)(0),$$
$$w_1 \otimes w_2 \mapsto o^{(i)}_{\mathcal Y}(w_1)w_2,$$
is an $A(L(c,0))$--module homomorphism.
\end{lemma}
\begin{proof}
Notice that both ${\mathcal Y}^{(0)}$ and ${\mathcal Y}^{(1)}$ satisfy 
the Jacobi identity (but the $L(-1)$--property does not hold for ${\mathcal
Y}^{(0)}$). To prove that $\pi^i({\mathcal Y})$ gives an $A(L(c,0))$--homomorphism
we can use the same proof as in \cite{FZ} (Lemma 1.5.2), since it uses
only the Jacobi
identity. Then
$$o(a)o_{\mathcal Y}^{(i)}(v)=o_{\mathcal Y}^{(i)}(a*v)$$ for $i=0,1.$ 
\end{proof}

Notice that there is no non--trivial intertwining operator such
that ${\mathcal Y}^{(0)} =0$. Otherwise 
\bea
&& {\mathcal Y}^{(1)}(L(-1)w, x){\rm log}(x)=\nn
&& =\frac{d}{dx} \left( {\mathcal Y}^{(1)}(w, x){\rm log}(x) \right) \nn
&& =(\frac{d}{dx}{\mathcal Y}^{(1)}(w, x)){\rm log}(x)+
\frac{1}{x}{\mathcal Y}^{(1)}(w, x) \nonumber.
\eea
The left hand side does not contain non--logarithmic operators
hence \\ ${\mathcal Y}^{(1)}(w, x)=0$, for every $w$.
The following lemma and its proof are motivated by Lemma 2.10 in \cite{Li1}.
\begin{lemma} \label{l1}
The mapping 
$${\pi} : I \ \binom{W_2(c,h)}{L(c,h_1) \ L(c,h_2)} \rightarrow
{\rm Hom}(A(L(c,h_1),L(c,h_2)) ,W_2(c,h)(0)),$$
$${\mathcal Y} \mapsto \pi^0_{{\mathcal Y}},$$
is injective.
\end{lemma}
\begin{proof}
Assume that $\pi^0_{{\mathcal Y}}=0$ for some ${\mathcal Y}$. Then
\begin{equation} \label{zero}
\langle w'_3,{\mathcal Y}^{(0)}(w_1,x)w_2 \rangle=0,
\end{equation}
for every $w'_3 \in W_2(c,h)'(0)=M_2(c,h)(0)$, 
$w_1 \in L(c,h_1)$ and $w_2 \in
L(c,h_2)(0)$.
From the commutator formula it follows that
$$\langle w_3',{\mathcal Y}^{(0)}(w_1,x)L(-n)w_2 \rangle=0,$$
and therefore (\ref{zero}) holds for every $w_2 \in L(c,h)$.
Hence $<w'_3,{\mathcal Y}^{(0)}(w_1,x)w_2 \rangle>=0$ 
for every $w_1 \in L(c,h_1)$ and 
$w_2 \in L(c,h_2)$. \\
{\em Claim:} ${\mathcal Y}^{(0)}=0$. 
If not then  $W':=\{ (w_1)^{(0)}_n w_2 \ : \ w_i \in L(c,h_i), 
i=1,2, n \in \mathbb{C} \}$ 
is a submodule of $W_2(c,h)$ such that $W'(0) \cap W_2(c,h)(0)=0$.
But, since $W_2(c,h)$ does not have singular  
vectors out side $W_2(c,h)(0)$ the same thing holds for $W'$.
Therefore $W'$ has to be zero.
Therefore ${\mathcal Y}^{(0)}=0$, which implies ${\mathcal
Y}=0$, by the above argument. 
\end{proof}

\begin{lemma} \label{l2}
$${\rm dim} \ I \ \binom{W_2(c,h)}{L(c,h_1) \ L(c,h_2)} \leq 2 .$$ 
\end{lemma}
\begin{proof}
From \cite{M1} it follows that 
$$A(L(c,h_1),L(c,h_2)) \cong \bigoplus_{i=1}^k \frac{\mathbb{C}[x]}{\langle (x-a_i)^{l_{a_i}}\rangle },$$
for some $k, l_i \in \mathbb{N}$ and $a_i \in \mathbb{C}$. On the
other hand $W_2(c,h)(0) \cong \frac{\mathbb{C}[x]}{\langle
(x-h)^2\rangle}$ where an isomorphism is given by $x-h \mapsto v, 1 \mapsto w$.
Hence if $l_{a_i} \geq 2$ the space of $A(L(c,0))$--homomorphisms
is at most two--dimensional.
By using the previous lemma we have the proof.
\end{proof}

Now the aim is to construct an intertwining operator starting from
$$\psi \in {\rm Hom}_{A(L(c,0)}(A(L(c,h_1) \otimes_{A(L(c,0))} L(c,h_2)(0),
W_2(c,h)(0)).$$
This can be exhibited (under some conditions) by using
Frenkel--Zhu--Li's construction (cf. \cite{FZ}, \cite{Li1}).
We stress that for our example, we use only a part
of the construction because we can use singular vectors.
Here is an outline (for more details see \cite{Li1}):
\begin{itemize} \label{long}
\item[(a)] To every vertex operator algebra $V$ one
associates a Lie algebra 
$g(V)=g_-(V) \oplus g_0(V) \oplus g_+(V)$. There is a
projection $g(V) \rightarrow A(V)$.
For every $g_0(V)$--module $U$ we consider a standard induced $g(V)$--module 
$F(U)$. In the case when $V=L(c,0)$ and $M=L(c,h)$, it is enough to
work with the Virasoro subalgebra of $g(L(c,0))$. Also 
$M(c,h)\hookrightarrow  F(L(c,h)(0))$, as a $Vir$--module
(which is a $L(c,0)$--module). 
Lift the homomorphism $\psi$ to a homomorphism (we keep the
same notation $\psi$) between $L(c,h_1) \otimes L(c,h_2)(0)$ and
$W_2(c,h)(0)$. 
\item[(b)] Define a $g$--module (where $g=g(L(c,0))$ )
structure on $$T :=\mathbb{C}[t,t^{-1}] \otimes L(c,h_1) \otimes
M(c,h) \hookrightarrow \mathbb{C}[t,t^{-1}] \otimes L(c,h_1) \otimes
F(L(c,h_2)(0))$$
and  a corresponding bilinear pairing between $M_2(c,h)(0)$ and $T$ by letting
$$\langle w'_3, t^n \otimes w_1 \otimes w_2 \rangle=\delta_{n-{\rm
deg}(w_1)+1,0}\langle w'_3,\psi (w_1 \otimes w_2)\rangle,$$
for $w_1 \in L(c,h_1)$ and $w_2 \in L(c,h_2)(0)$. Extend this pairing
for every $w_2 \in L(c,h_2)$.
Similarly we can extend this pairing to $M_2(c,h)$ and $T$ such that
for every $y \in U({\rm Vir})$
$$\langle yf,v\rangle=\langle f,\theta(y)v\rangle,$$
where $\theta$ is an anti--involution of $U({\rm Vir})$ (in particular
$\theta(L(n))=L(-n)$).
Set $${\mathcal Y}_t(w_1,x)=x^{h_1-h_2-h_3} \sum_{n \in \mathbb{Z}} (t^n
\otimes w_1 )x^{-n-1}.$$
\item[(c)] Prove the commutativity and associativity for ${\mathcal
Y}_t$ (this is the hardest part of the construction).
\item[(d)] Show that $L(-1)$--property holds for ${\mathcal Y}_t$ (this
requires that $W_3$ is 
$L(0)$--diagonalizable). 
\item[(e)] Define 
$${\mathcal Y}(w_1,x)w_2:=\psi({\mathcal Y}_t(w_1,x) \otimes w_2),$$
for which $L(-1)$--property might not holds, of the type 
\begin{equation} \label{zzz}
\binom{M_2(c,h)'}{L(c,h_1) \ M(c,h_2)},
\end{equation}
(note that  $M_2(c,h) \hookrightarrow F(M_2(c,h)(0))$). 
By using an argument as in \cite{M1} it follows that
every intertwining operator of the type (\ref{zzz}) can be pushed down to an
intertwining operator of the form $$\binom{W_2(c,h)}{L(c,h_1)
\ M(c,h_2)}.$$
\item[(f)] Show that $L(-1)$--property holds.
\item[(g)] Check whether the intertwining operator projects to an intertwining
operator of the type
$$\binom{W_2(c,h)}{L(c,h_1)
\ L(c,h_2)}.$$
\end{itemize}

Suppose that
$$[A(L(c,h_1),L(c,h_2)):W_2(c,h)(0)]=1,$$
where $[M_1:M_2]$ stands for the multiplicity of $M_2$ inside $M_1$.
Therefore 
$${\rm Hom}_{A(L(c,0))}(A(L(c,h_1) \otimes_{A(L(c,0)} L(c,h_2)(0),
W_2(c,h)(0))$$
is two dimensional. Let us pick a surjective homomorphism $\psi^{(0)}$ (such
exists!) and the homomorphism $\psi^{(1)}=(L(0)-h)\psi^{(0)}$ (it has
one--dimensional image). If we apply procedure (a)-(e) for these homomorphisms
we obtain ``intertwining operators'' \footnote{
However, ${\mathcal Y}^{(0)}$ fails to satisfy $L(-1)$--property.} 
${\mathcal Y}^{(0)}$ and ${\mathcal Y}^{(1)}$, respectively. 
On the other hand ${\mathcal Y}^{(1)}$ is a genuine intertwining operator
of the type $\binom{L(c,h)}{L(c,h_1)
\ M(c,h_2)},$ via $L(c,h) \hookrightarrow W_2(c,h)$.
\begin{lemma} \label{L(-1)}
If ${\mathcal Y}^{(0)}$ and ${\mathcal Y}^{(1)}$ are as the above
$${\mathcal Y}^{(0)}(L(-1)w_1,x)=\frac{d}{dx}{\mathcal
Y}^{(0)}(w_1,x)+\frac{1}{x}{\mathcal Y}^{(1)}(w_1,x),$$
for every $w_1 \in L(c,h_1)$.
\end{lemma}
\begin{proof}
We prove first that 
\bea
&& \langle w'_3,{\mathcal Y}^{(0)}(L(-1)w_1,x)w_2\rangle= \nn
&& =\frac{d}{dx}\langle w'_3,{\mathcal Y}^{(0)}(w_1,x)w_2\rangle+\frac{1}{x}\langle w'_3,{\mathcal
Y}^{(1)}(w_1,x)w_2\rangle,
\eea
for every $w_1 \in L(c,h_1), w_3 \in W_2(c,h)(0)$ and $w_2 \in M(c,h_2)$.
Note that 
$$\langle w'_3,(L(0)-h_3){\mathcal Y}^{(0)}(w_1,x)w_2\rangle=\langle w'_3,{\mathcal Y}^{(1)}(w_1,x)w_2\rangle.$$
Then by the Jacobi identity
\bea
&& \langle w'_3,{\mathcal Y}^{(0)}(L(-1)w_1,x)w_2\rangle= \\
&& =\frac{1}{x}\langle w'_3,L(0){\mathcal
Y}^{(0)}(w_1,x)w_2\rangle \nn
&& -\frac{{\rm wt}(w_1)+
{\rm wt}(w_2)}{x} \langle w'_3,{\mathcal Y}^{(0)}(w_1,x)w_2\rangle \nn
&& = \frac{1}{x}\langle w'_3,{\mathcal Y}^{(1)}(w_1,x)w_2\rangle \nn
&& + \frac{h_3-{\rm wt}(w_1)-{\rm wt}(w_2)}{x}\langle w'_3,{\mathcal
Y}^{(0)}(w_1,x)w_2\rangle \nn
&& = \frac{1}{x}\langle w'_3,{\mathcal Y}^{(1)}(w_1,x)w_2\rangle+\frac{d}{dx}
\langle w'_3,{\mathcal Y}^{(0)}(w_1,x)w_2\rangle.\nonumber
\eea
As in \cite{Li1}, Lemma 3.9, we check that for every $m \in
\mathbb{N}$ and $w'_3 \in W'_3(0)$
\bea \label{L(-1)1}
&& \langle w'_3,L(m){\mathcal Y}^{(0)}(L(-1)w_1,x)w_2\rangle = \nn
&& = \frac{d}{dx}\langle w'_3,L(m){\mathcal Y}^{(0)}(w_1,x)w_2\rangle+
\frac{1}{x}\langle w'_3,L(m){\mathcal Y}^{(1)}(w_1,x)w_2\rangle.
\eea
>From this fact it follows that
\bea \label{L(-1)2}
&& \langle w'_3,{\mathcal Y}^{(0)}(L(-1)w_1,x)w_2\rangle= \nn
&& = \frac{d}{dx}\langle w'_3,{\mathcal Y}^{(0)}(w_1,x)w_2\rangle+
\frac{1}{x}\langle w'_3,{\mathcal Y}^{(1)}(w_1,x)w_2\rangle,
\eea
for every $w'_3 \in W_3$. So we have the proof.
\end{proof}
\begin{theorem} \label{main}
Suppose that 
\begin{equation} \label{main1}
[A(L(c,h_1),L(c,h_2)):W_2(c,h)(0)]=1.
\end{equation}
Then $I \binom{W_2(c,h)}{L(c,h_1) M(c,h_2)}$ is two dimensional. 
Moreover, if 
\begin{equation} \label{main2}
[A(L(c,h_2),L(c,h_1)):W_2(c,h)(0)]=1
\end{equation}
as well, then $\pi$ in Lemma \ref{l1} is an isomorphism and
$I \binom{W_2(c,h)}{L(c,h_1) \ L(c,h_2)}$ is two dimensional. 
In particular, there is a non--trivial logarithmic intertwining
operator of type $\binom{W_2(c,h)}{L(c,h_1) \ L(c,h_2)}$.
\end{theorem}
\begin{proof}
If (\ref{main1}) holds then we have constructed ${\mathcal Y}^{(0)}, {\mathcal
Y}^{(1)} \in I \ \binom{W_2(c,h)}{L(c,h_1) \ M(c,h_2)}$, such that
the $L(-1)$--property does not hold for ${\mathcal Y}^{(1)}$. 
If we let 
$${\mathcal Y}(w_1,x)={\mathcal Y}^{(0)}(w_1,x)+{\mathcal
Y}^{(1)}(w_1,x){\rm log}(x),$$
then ${\mathcal Y}$ clearly satisfies the Jacobi identity. From 
Lemma \ref{L(-1)} it follows that 
$${\mathcal Y}(L(-1)w_1,x)=\frac{d}{dx}{\mathcal Y}(w_1,x).$$
Hence ${\mathcal Y}$ is a logarithmic intertwining operator
of type $\binom{W_2(c,h)}{L(c,h_1) \  M(c,h_2)}.$
Now in order to project ${\mathcal Y} \in I \ \binom{W_2(c,h)}{L(c,h_1) \ M(c,h_2)}$
to $\bar{{\mathcal Y}} \in I \ \binom{W_2(c,h)}{L(c,h_1) \ L(c,h_2)}$
one has to check that
$$\langle w'_3, {\mathcal Y}(w_1,x)w_2\rangle=0,$$
for every $w_2 \in M'(c,h_2)$, where $M'(c,h_2)$ is the
maximal submodule of $M(c,h_2)$. 
Since $M'(c,h_2)$ is generated by the singular
vectors and the Jacobi identity holds, it is enough to check that
$$\langle w'_3, {\mathcal Y}(w_1,x)v_{sing}\rangle=0,$$
for $w'_3 \in W_2(c,h)$ and $w_1 \in L(c,h_1)(0)$. 
Write $v_{sing}=Pv_{h_2}$ for some
$P \in U(Vir)_-$. Then
$$\langle w'_3, {\mathcal Y}(w_1,x)Pv_{h}\rangle=R(\partial_x)\langle w'_3, {\mathcal
Y}(w_1,x)v_{h}\rangle,$$
for some differential operator $R(\partial_x)$ (cf. Section 1.5).
This Euler's differential equation reduces to
a differential equation with constant coefficients.
Because of the condition (\ref{main2})
the corresponding characteristic equation has double roots
at $h$, there is a logarithmic solution.
Therefore
$$P(\partial_x)x^{h-h_1-h_2}=R(\partial_x){\rm
log}(x)x^{h-h_1-h_2}=0,$$
and hence $\langle w'_3, {\mathcal Y}(w_1,x)v_{sing}\rangle=0$, for every $w'_3 \in
W_2(c,h)$.
\end{proof}

\begin{remark}
Notice that in general
$$A(L(c,h_1),L(c,h_2)) \ncong A(L(c,h_2),L(c,h_1)).$$

\end{remark}

\section{Applications}

Here we illustrate some consequences of Theorem \ref{main}.
When $c=1$ we show that we can obtain certain logarithmic 
intertwining operators but only in the case when one of
the modules is a Verma module, i.e. the second condition
in Theorem \ref{main} does not hold. 
In the case when $c=-2$ the
situation is substantially different where, for a specific model,
both conditions hold.

\subsection{$c=1$}

This model was closely examined in \cite{M1} where we calculated
$$A(L(1,\frac{m^2}{4}),L(1,\frac{n^2}{4}))$$ for every pair $m,n \in
\mathbb{N}$. 
The result is:
\begin{lemma}[M]
$$A(L(1,\frac{m^2}{4}),L(1,\frac{n^2}{4})) \cong 
\frac{\mathbb{C}[x]}{\langle \prod_{i \in J_{m,n}} (x-\frac{i^2}{4})\rangle}.$$
where $J_{m,n}$ is the multi set $\{m+n,m+n-2,...,m-n \}$.
\end{lemma}
\begin{corollary} \label{c=1}
Suppose that $m > n$. Then for every $i$, satisfying $0< i < |n-m|$, there is
a logarithmic intertwining operator of type
$$\binom{W_2(1,\frac{i^2}{4})}{L(1,\frac{m^2}{4}) \ 
M(1,\frac{n^2}{4})}.$$
Consequently there are {\em no} nontrivial logarithmic operators of
type 
\begin{equation} \label{jedan}
\binom{W_2(1,\frac{i^2}{4})}{L(1,\frac{m^2}{4})
\ L(1,\frac{n^2}{4})},
\end{equation}
i.e., the space of logarithmic intertwining operators
of type (\ref{jedan}) is one dimensional.
\end{corollary} 
\begin{proof}
Directly follows from Theorem \ref{main} and Lemma \ref{c=1}.
\end{proof}

\subsection{$c=-2$}

It was shown
in \cite{Gu1} that certain matrix coefficients built up from primary fields of
the conformal weight $h=\frac{-1}{8}$ with $c=-2$ 
yield logarithmic singularities.
We consider
a vertex operator algebra $L(-2,0)$. 
\begin{lemma} \label{-2}
$$A(L(-2,\frac{-1}{8}), L(-2,\frac{-1}{8})) \cong
\frac{\mathbb{C}[x]}{\langle x^2\rangle}.$$
\end{lemma}
\begin{proof}
It is not hard to see that
$$A(L(-2,\frac{-1}{8})) \cong \frac{\mathbb{C}[x,y]}{\langle p(x,y)\rangle},$$
where $\langle p(x,y)\rangle$ is a $\mathbb{C}[y]$--(bi)module generated by
$$p(x,y)=y^2-2xy-\frac{y}{4}+x^2-\frac{x}{4}-\frac{3}{64}.$$
%
%
Hence 
$$A(L(-2,\frac{-1}{8}), L(-2,\frac{-1}{8})) \cong
\frac{\mathbb{C}[x]}{\langle x^2\rangle}.$$
\end{proof}

\begin{corollary} \label{c=-2}
The space 
$$I \binom{ W_2(-2,0)}{L(-2,\frac{-1}{8})
\ L(-2,\frac{-1}{8})},$$
is two dimensional.
\end{corollary}
\begin{proof}
Follows directly from Lemma \ref{-2} and Theorem (\ref{main}).
\end{proof}

\subsection{ Why is $c=0$ different ?}

\begin{figure}
\begin{center}
\begin{picture}(110,285)
\put(-12,180){$\scriptstyle{M(0,1)}$}
\put(88,180){$\scriptstyle{M(0,2)}$}
\put(-12,275){$\scriptstyle{M(0,0)}$}

\multiput(0,0)(0,14.5){6}{\line(0,1){12}}
\multiput(100,0)(0,14.5){6}{\line(0,1){12}}

\multiput(90,0)(-12,12){7}{\line(-1,1){10}}
\multiput(10,0)(12,12){7}{\line(1,1){10}}

\put(0,95){\vector(0,1){80}}
\put(100,95){\vector(0,1){80}}

\put(10,95){\vector(1,1){80}}
\put(90,95){\vector(-1,1){80}}

\put(0,190){\vector(0,1){80}}
\put(90,190){\vector(-1,1){80}}
\end{picture}
\end{center}
\caption{\label{em00}Embedding structure for $M(0,0)$} 
\end{figure}
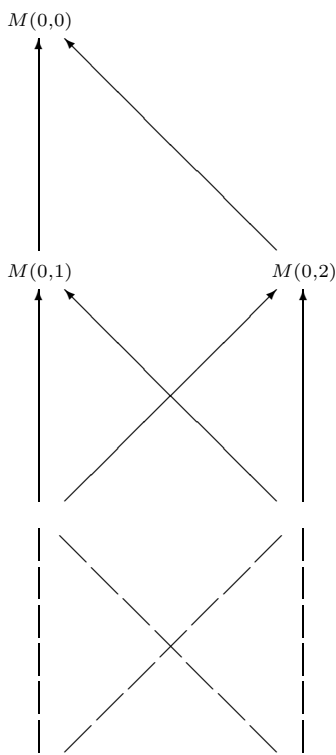

Note that $c_{2,3}=0$, and thus we are dealing with the minimal models.
$L(0,0)$ is the trivial vertex operator algebra. Therefore the only
irreducible module is $L(0,0)$ itself. To avoid these trivialities
we consider a vertex operator algebra $M_0=M(0,0)/\langle L(-1){v_0}\rangle$
(which is not simple !).
The embedding structure for $M(0,0)$ is given by (\ref{em00}) (cf. \cite{RW}).
Set 
$$h_m=\frac{m^2-1}{24}.$$ Then it is known (see \cite{RW}) that
$M(0,h_m)$ is reducible. For computational purposes we will
study more closely irreducible $M_0$--modules associated to an infinite
sequence $m=3p-2$ and the corresponding $M_0$--modules $L(0,h_m)$ (all results
can be easily derived for every $L(0,\frac{m^2-1}{24}$)). 

If $p=2$ and $q=3$ then
one sees (after some calculations) that
\begin{equation} \label{32}
A(L(0,\frac{5}{8}),L(0,\frac{5}{8})) \cong
\frac{\mathbb{C}[x]}{\langle x(x-2)\rangle}.
\end{equation}
More generally, for every even $p$ we have
\bea \label{03p}
&& A(L(0,\frac{(3p-2)^2-1}{24}),L(0,\frac{(3p-2)^2-1}{24}))
\cong \nn
&& \frac{\mathbb{C}[x]}{
\prod_{\frac{1}{2} \leq n \leq \frac{p-1}{2}}
\left(\frac{(3p-1-6n)
(3p-3-6n)}{24}-x\right)\left(\frac{(3p-1+6n)(3q-3+6n)}{24}-x \right)}.
\eea
If we analyze (\ref{03p}) \footnote{Notice that all numbers that
appear in the polynomial in (\ref{03p}) are pentagonal
numbers of the form $\frac{3k^2+k}{2}$.}
more carefully we see that 
the right hand side is divisible by $x$ but not by $x^2$ (compare with the case
$c=-2$). Hence our Theorem \ref{main} {\rm does not} yield 
any non--trivial logarithmic 
intertwining operators. This is because only the top level 
is relevant in our approach.
In other words if the top level is one dimensional,  
Theorem \ref{main} does not give a sufficient condition
for the existence of a logarithmic
intertwining operators. Logarithmic operators are invisible.
Also notice that $W_2(0,0)$ has a degenerate structure and
it corresponds to a semi--direct product of  a trivial module
and $M(0,0)$, i.e. we have
$$0 \rightarrow \mathbb{C} \rightarrow W_2(0,0) \rightarrow M(0,0)
\rightarrow 0.$$ 

Still in  \cite{GL}  some interesting logarithmic conformal field
theories were considered exactly for the models we discussed 
above.  

To understand their result 
let us consider a small detour into conformal field theory.
If $A(z)$ is a primary field and of conformal weight $h$ and $c \neq
0$ then physicists often write the following ``operator product expansion''(OPE):
\begin{equation} \label{5}
A(z)A(0)=z^{-2h}(1+\frac{2h}{c}T(0)z^2+ \cdots ),
\end{equation}
where $T(z)$ is the Virasoro field.
The formula (\ref{5}) can be easily reformulated:
\begin{proposition} \label{exp5}
Suppose that ${\mathcal Y} \in I \binom{L(c,0)}{L(c,h) \  L(c,h)}$, $c
\neq 0$ and ${\mathcal Y}$ is normalized such that
$${\mathcal Y}(w_{1},x)w_{1}=x^{-2h}{\bf 1}+ax^{-2h+2}L(-2){\bf
1}+\ldots,$$ 
where $w_1$ is the lowest weight vector of $L(c,h)$.
Then $a=\frac{2h}{c}$.
\end{proposition}
\begin{proof}
Directly follows from the formulas:
$$\langle L(-2){\bf 1}',{\mathcal Y}(w_{1},x)w_{2}\rangle=a\frac{c}{2}$$
and 
$$[L(2),{\mathcal Y}(w_1,x)]=(x^3 \frac{d}{dx}+3hx^2){\mathcal Y}(w_1,x).$$
\end{proof}

It was observed in \cite{GL} and \cite{Gu2}
that one can reformulate (\ref{5}), by adding certain logarithmic operators,
such that it makes sense even when $c=0$.
Then the following
OPE was proposed (from now on $c=0$):
\begin{equation} \label{6}
A(z)A(0)=z^{-2h}(1+\frac{2h}{b}z^{-2h}\left(t(0)+{\rm log}(z)T(0)+...\right),
\end{equation}
for some constant $b$. Note that the first logarithmic operator
appears on the weight $2$. Note that the Virasoro vector is contained
in $V(2)$ as well. As noticed in \cite{GL},
this makes the whole $c=0$ theory peculiar.
OPE (\ref{6}) implies that there should exists a logarithmic intertwining operator
\begin{equation} \label{lastfor}
{\mathcal Y} \in I \ \binom{W_b}{L(c,\frac{m^2-1}{24}) \ L(c,\frac{m^2-1}{24}) },
\end{equation}
for some logarithmic module $W_b$, that depends on some constant $b$, such that
\begin{equation} \label{2vacuums}
{\mathcal Y}(w_{1},x)w_{1}=x^{-2h}{\bf 1}+\frac{2h}{b}({\rm
log}(x)x^{-2h+2}L(-2){\bf 1}+tx^{-2h+2}+ ...,
\end{equation}
where $t  \in W_b$ satisfies $L(2)t=\frac{b}{2} {\bf 1}$.
It is not clear what kind of logarithmic module is
$W_b$.  Let us try to construct such a module starting from $M_2(0,0)$. 
(\ref{2vacuums}) predicts that $W_b(0)$ is one--dimensional, $W_b(1)$
is trivial and $W_b(2)$ is spanned by
$\omega=L(-2)v$ and some vector $t$.
Hence the structure of $W_b$ is given by
Figure (\ref{m00})
(arrows indicate where corresponding operators map one--dimensional spaces).
\begin{figure}
\begin{center}
\begin{picture}(102,195)
\put(-12,87){$\scriptstyle{L(-2)v}$}
\put(98,87){$\scriptstyle{t}$}
\put(-2,186){$\scriptstyle{v}$}
\put(-24,140){$\scriptstyle{L(-2)}$}
\put(50,145){$\scriptstyle{L(2)}$}

\put(0,0){\line(0,1){80}}
\put(100,0){\line(0,1){80}}

\put(0,180){\vector(0,-1){80}}
\put(92,100){\vector(-1,1){80}}
\end{picture}
\caption{\label{m00} } 
\end{center}
\end{figure}
If we analyze possible choices we realize that one cannot construct
such subquotient. Also one cannot accomplish that 
$W_b(1)=0$.

\subsection{Towards $c=0$ logarithmic intertwiners}

In the previous section we argue that one cannot 
recover $W_b$ as a subquotient of $M(0,0)$. 
This is because our considerations were solely on the level of the Virasoro 
algebra. We shall show that one can 
consider a larger Lie algebra and a module which exhibit certain
properties similar to---yet to be constructed---$W_b$.
Let
$$W_{{\rm log}} = \bigoplus_{i \in \mathbb{Z}, m \in \mathbb{Z}} \mathbb{C}t^{(i)}(m).$$ 
Then it is not hard to see that the relation
\begin{equation} \label{logvir}
[t^{(i)}(m),t^{(j)}(n)]=(m-n)t^{(i+j)}(m+n)+(j-i)t^{(i+j-1)}(m+n),
\end{equation}
where $i, j \in \mathbb{Z}$ and $m,n \in \mathbb{Z}$,
defines a Lie algebra structure on $W_{{\rm log}}$ (``logarithmic algebra'').
The Lie algebra $W_{{\rm log}}$ has the following representation
in terms of formal logarithmic vector fields:
$$t^{(i)}(n) \mapsto x^{-n+1}{\rm log}^i(x)\frac{d}{dx}.$$
As usual we are interested in central extensions.
\subsection{Central extension of $W_{{\rm log}}$}
We will prove the following theorem.
\begin{theorem}
$${\rm dim} \ H^2(W_{{\rm log}},\mathbb{C})=1.$$
Moreover a non--trivial 2--cocycle $c( \ , \ )$ is given by
$$ c(t^{(1-i)}(-m),t^{(1-j)}(-n))=
\sum_{r \geq 0}^{i+j} \frac{m^r n^{i+j-r} ((i-r)^3-(i-r))}{12(j+i-r)!
r!}.$$
\end{theorem}

\begin{lemma} \label{firstl}
Suppose that ${\mathcal L} \subset {\mathbb C}((t))\frac{d}{dt}$
where ${\mathbb C}((t))\frac{d}{dt}$ is equipped with a
Lie algebra structure in the natural way. 
Suppose that the following holds:
\begin{itemize} 
\item[$(a)$] 
$$t^n \frac{d}{dt} \in {\mathcal L}$$
for every $n \in \mathbb{Z}$, i.e.,
$W$ contains a (polynomial) Virasoro subalgebra.
\item[$(b)$]
If $f(t) \frac{d}{dt} \in {\mathcal L}$, 
$f(t)=\sum_{i} a_i t^i \in \mathbb{C}((t))$, and
${\rm Res}_t f(t)=0$, then 
$$\int f(t) dt:=\sum_i a_i \frac{t^{i+1}}{i+1} \in {\mathcal L}.$$
\end{itemize}
Then 
$${\rm dim} \ H^2({\mathcal L},\mathbb{C})=1 \ {\rm or} \ 0.$$
\end{lemma}

\begin{proof}
It follows from the proof of the main
result in \cite{LW}.  This proof was used
to show that $\mathbb{C}((t))\frac{d}{dt}$
has essentially unique central extension.
\end{proof}
The previous lemma implies that it is enough to
construct a non--trivial 2--cocycle on ${\mathcal L}$ to prove
that $H^2({\mathcal L}, \mathbb{C})$ is non--trivial.
\begin{lemma}
A Lie subalgebra ${\mathcal S} \subset \mathbb{C}((t)) \frac{d}{dt}$
spanned by $t^{-i+1} e^{mt} \frac{d}{dt}$, where $i, m \in \mathbb{Z}$.
Then ${\mathcal S}$ satisfies conditions in Lemma \ref{firstl}.
Moreover, ${\mathcal S} \cong W_{\rm log}$.
\end{lemma}
\begin{proof}
It is trivial to prove that ${\mathcal S}$ is a Lie subalgebra.
Now, integration by parts implies ``integral'' condition (b).
If we take $m=0$ we obtain the Virasoro subalgebra.
If we let $\tilde{t}^{(i)}(m)=t^{-i+1}e^{mt} \frac{d}{dt}$
then we have 
$$[\tilde{t}^{(i)}(m).\tilde{t}^{(j)}(n)]=(i-j)\tilde{t}^{(i+j)}
(m+n)+(n-m)\tilde{t}^{(i+j-1)}(m+n).$$
Then $$\tilde{t}^{(i)}(m) \mapsto t^{(1-i)}(-m),$$
is a Lie algebra isomorphism.
\end{proof}
\begin{proof}
Let $\tilde{c}({\tilde t} ^{(i)}(m),{\tilde t}^{(j)}(n)),$
be a 2--cocycle considered in \cite{LW}
(where was denoted by $\alpha$--which is essentially 
the Gelfand--Fuchs cocycle)
and $c( \ , \ )$ an non--trivial 2--cocycle on 
$W_{\rm log}$ induced by the isomorphism $\varphi$. 
Lemma \ref{firstl}
implies that ${\rm dim} \ H^2({\mathcal S},\mathbb{C})=1.$
Hence  ${\rm dim} \ H^2(W_{{\rm log}},\mathbb{C})=1.$
\end{proof}

Here are some properties of $W_{{\rm log}}$.
\begin{proposition}
\begin{itemize}
\item[$(a)$]
${W}_{{\rm log}}$ is generated by $t^{(-1)}(m)$, $t^{(0)}(n)$ and 
$t^{(1)}(p)$,  where $m,n,p  \in \mathbb{Z}$.
\item[$(b)$]
The Lie algebra ${W}_{{\rm log}}$ has a triangular--like decomposition
$${W}_{{\rm log}}={W}^+_{{\rm log}} \oplus {W}^0_{{\rm log}} \oplus
{W}^-_{{\rm log}},$$
where 
$${W}^+_{{\rm log}}=\bigoplus_{i >0, m \in \mathbb{Z}} \mathbb{C} t^{(i)}(m),$$
$${W}^-_{{\rm log}}=\bigoplus_{i <0, m \in \mathbb{Z}} \mathbb{C} t^{(i)}(m),$$
$${W}^0_{{\rm log}}=\bigoplus_{i \in \mathbb{Z}} \mathbb{C} t^{(0)}(i) \oplus \mathbb{C}b.$$
\item[$(c)$]
${W}_{{\rm log}}$ has an anti--involution defined by 
$$t^{(i)}(m) \mapsto (-1)^i t^{(i)}(-m).$$
\end{itemize}
\end{proposition}
\begin{proof}
The proof is straightforward.
\end{proof}
Let us denote by $\hat{W}_{{\rm log}}=W_{{\rm log}} \oplus \mathbb{C}b$ 
the Lie algebra
induced by 2--cocycle $c( \ , \ )$. We denote the generator of
one--dimensional central subspace by $b$.
If $i=j=0$, then (\ref{logvir}) gives us commutation relation
for the Virasoro algebra with the {\em trivial} central charge.
This is in the spirit of \cite{GL}.
It is not true that central charge is absent. Instead it is
trivial for the Witt subalgebras.

We consider a Verma module
$$V(b,0):={\rm Ind}^{\hat{W}_{{\rm log}}}_{{W}^+_{{\rm log}} \oplus
{W}^0_{{\rm log}} \oplus \mathbb{C}b} \mathbb{C}{v_b},$$
where ${\mathbb C}v_b$ is an one--dimensional module such that 
$b.v_b=bv_b$ (on the right hand side $b \in \mathbb{C}$)
and the action for remaining generators is trivial.
Notice that this module is well defined because 
$$[{W}^+_{{\rm log}} \oplus {W}^0_{{\rm log}},{W}^+_{{\rm log}} \oplus {W}^0_{{\rm log}}]|_{\mathbb{C}{v_b}}=0.$$

Now our discussion is heuristic.
Let $W(b,0)$ be some subquotient of the Verma module $V(b,0)$. 
$t:=t^{(-1)}(-2){v_b}$  and $\omega=t^{(0)}(-2){v_b}$ 
are nontrivial. Then 
$t^{(-1)}(-2){v_b}$ is a logarithmic operator but {\em not}
for the horizontal Virasoro generator $t^{(0)}(0)$, but rather
for the vertical Virasoro generator $t^{(1)}(0)$.  
Now let us go back to our construction. 
If ${\mathcal Y}$ is such an operator (normalized such that
${\mathcal Y}(w_1,x)w_2={\bf 1}x^{-2h}+\cdots $)
then the (formal) matrix coefficient:
\begin{equation} \label{t2}
\langle t^{(-1)}(-2){v_b}', {\mathcal Y}(w_1,x)w_2\rangle
\end{equation}
should include some logarithmic terms.
As before we can determine (\ref{t2}) by
solving an appropriate differential equation. This time
the equation is
\begin{equation} \label{odeb}
P(\partial_x)\langle t^{(-1)}(-2){v_b}', {\mathcal
Y}(w_1,x)w_2\rangle=Q(\partial_x)x^{-2h}.
\end{equation}
\begin{example} \label{lastex}
Let us consider a special case: $L(0,\frac{5}{8})$.
Clearly, 
$$\langle t^{(-1)}(-2){v_b}',t^{(0)}(-2){v_b}\rangle=\frac{2b}{3}$$
and
$$\langle t^{(-1)}(-2){v_b}',t^{(1)}(-2){v_b}\rangle=0.$$
By analyzing singular vectors we obtain the following differential equation 
\bea \label{odec}
&& \left(\frac{d^2}{d^2 x}+\frac{3}{2}
x^{-1}\frac{d}{dx}-\frac{15}{16}x^{-2} \right)
\langle t^{(-1)}(-2){v_b}', {\mathcal Y}(w_1,x)w_2\rangle=\frac{2b}{3}x^{-\frac{5}{4}}.
\eea
Every (formal) solution of (\ref{odec}) is
of the form
\begin{equation} \label{eq}
\lambda x^{3/4}+\frac{b}{3}x^{3/4}{\rm log}(x), 
\end{equation}
By using the same argument as in Proposition \ref{exp5} 
we see that
$${\mathcal Y}(w_1,x)w_2=x^{-5/4}{\bf 1}+
x^{3/4}(\frac{2h}{b}t^{(1)}(-2){\bf 1}+\mu {\rm log}(x)t^{(0)}(-2){\bf 1})+\cdots$$
By combining this with (\ref{eq}) it follows that
$\mu=\frac{1}{2}$. Hence if we want (\ref{6}) to hold
then $\frac{2h}{b}=\frac{1}{2}$. Therefore $b=\frac{5}{2}$.
\end{example}

For all other modules $L(0,\frac{(3p-2)^2-1}{24})$, where $p$ is even,
and in general for every $L(0,\frac{m^2-1}{24})$ where
$$m \cong 0,2,3,4 \ {\rm mod} \ 6,$$
the situation is as follows:
There is a singular vector $v_{sing}$ 
inside $M(0, \frac{m^2-1}{24})$ that generates its  maximal submodule 
$M'(0,\frac{m^2-1}{24})$.
Hence, if we write
$$F(x)=\langle t^{(-1)}(-2){\bf 1}', {\mathcal Y}(w_1,x)w_2 \rangle,$$
the equation (\ref{odeb}) can be written as
\begin{equation} \label{lal}
P(\partial_x)F(x)=Q(\partial_x)x^{-2h},
\end{equation}
for some differential operators $P$ and $Q$.
By analyzing (\ref{03p}) and related formulas, 
we see that (\ref{lal}) has a logarithmic
solution if and only if $Q(\partial_x)x^{-2h} \neq 0$.
Still the expression for
$Q(\partial_x)x^{-2h}$ 
is unknown in general.
Some computations show (cf. \cite{GL}) that
the equality (\ref{6}) can be accomplish only for 
some special $b$ (we stress that our value $b$ differs
from the one in \cite{GL}).

\setcounter{equation}{0}
\section{Final remarks and some open problems}

\begin{itemize}
\item[(a)]  Determine complete embedding structure for $M_n(c,h)$ and
calculate graded characters for $\tilde{W}_n(c,h)$. 
\item[(b)] We know that
$M_0=M(0,0)/\langle L(-1){\bf 1}\rangle$ is a vertex operator algebra.
What algebraic structure, besides being a $M_0$--module, governs 
\begin{equation} \label{logfi}
M(0,0)/\langle L(-2){\bf 1}-\frac{3}{2}L(-1)^2 {\bf 1}\rangle \  ?
\end{equation}
\item[(c)]  Construct a non--trivial intertwining operators of the form
$$\binom{W_b}{L(\frac{m^2-1}{24},0)  \ L(\frac{m^2-1}{24},0)}.$$
We believe that this construction can be obtained by
constructing a refinement of Frenkel-Zhu's formula which 
carries an additional information about the module (not only the
top level). We shall study this problem in a sequel.

\end{itemize}

\renewcommand{\theequation}{\thesection.\arabic{equation}}
\setcounter{equation}{0}

\section{Appendix: Shapovalov form}

In this appendix we present some calculations of singular vectors in
the case of $M_2(c,h)$.
Let us define Shapovalov form $( \ , \ )$ on $M_2(c,h) \times M_2(c,h)$ in the
following way:
Fix $v,w  \in M_2(c,h)(0)$ as before
such that  $(v,v)=1$ and $(w,w)=1$ and $(v,w)=(w,v)=0$.
Then for every $s_1,s_2 \in U(Vir_-)$ we define a bilinear (it is not
symmetric) form on $M_2(c,h)$
\begin{equation} \label{shapform}
(s_1.u_1,s_2.u_2):=(u_1, s_1^T s_2 .u_2),
\end{equation}
where $u_1$ and $u_2 \in M_2(c,h)(0)$  
and $s^T= L(i_k) \ldots L(i_1)$ for $s=L(-i_1)\ldots L(-i_k)$.
Note that this form is not invariant in general.
Still we have the following result:
\begin{proposition}
The (right) radical
$${\rm Rad}(M_2(c,h)):=\{w_2: \ (w_1,w_2)=0 \ {\rm for} \ {\rm every}
\ w_1 \in M_2(c,h) \},$$
of the Shapovalov form (\ref{shapform})
is equal to  $M'_2(c,h) \subset M_2(c,h)$. 
\end{proposition}

As in the case of Shapovalov form 
for Verma modules (cf. \cite{FF1}, \cite{FF2})
we introduce matrices $[S_2(c,h)_n]$, $n \in \mathbb{N}$.
and the corresponding determinants  $S_2(c,h)_n:={\rm det}[S_2(c,h)_n]$. 
$[S_2(c,h)_n]$ is a matrix of size $2p(n)$ ($p(n)$ is
the number of partitions of $n$). 
$[S_2(c,h)_n]$ is a $2 \times 2$ block matrix with four
blocks of size $p(n) \times p(n)$, where the diagonal blocks are
equal to $S(c,h)_n$ (this is a matrix in the case
of ordinary Verma module $M(c,h)$) and one of the blocks
is a zero matrix (we have chosen a basis that involves the first $p(n)$ vectors from
$U(Vir)^n_-v$ and then the $p(n)$ vectors in $U(Vir)_-^n w$).
Hence $$S_2(c,h)_n=(S(c,h)_n)^2.$$
We shall need the following important lemma that describes the
remaining block in $[S_2(c,h)_n]$:
Let us assume that the Shapovalov form is $\mathbb{C}[h]$--valued, i.e. 
$h$ is a formal variable. Then we have:
\begin{lemma} \label{stupid}
For every $s_1, s_2 \in U(Vir)_-$,
$$(s_1.v,s_2.w)=\frac{d}{dh}(s_1.v,s_2.v).$$
In other words $[S_2(c,h)_n]$ can be written as a block matrix
$$\left(\begin{array}{cc}[S(c,h)] & \frac{d}{dh} [S(c,h)] \\ 0 & [S(c,h)]
\end{array}\right).$$
\end{lemma}
\begin{proof}
denote by ${\rm \pi}$ the projection to a subspace 
$U(Vir)_-.v$.
It is enough to show that
\begin{equation} \label{induction}
{\pi} L(j)L(-i_1)\cdots L(-i_k)w=\frac{d}{dh}L(j)L(-i_1)\cdots L(-i_k)v,
\end{equation}
for every $j,k \geq 0$, $i_1,...,i_k \in \mathbb{N}$.
If $j=0$ the result holds. Suppose that (\ref{induction}) holds
for $j=0,...,m$. Then 
\bea 
&&{\pi} L(m+1)L(-i_1)\cdots L(-i_k)w=\\
&& ={\pi}( (m+1+i_1)L(m+1-i_1)L(-i_2)\cdots
L(-i_k)w \nn
&& +(m+1+i_2)L(-i_1)L(m+1-i_2)\ldots L(-i_k)w  \nn
&& +...+ (m+1-i_k)L(-i_1)\ldots L(-i_{k-1})L(m+1-i_k)w) \nn
&& =\frac{d}{dh}(m+1+i_1)L(m+1-i_1)L(-i_2)\ldots L(-i_k)v \nn
&& +...+ (m+1-i_k)\frac{d}{dh} L(-i_1)\ldots L(-i_{k-1})L(m+1-i_k)v \nn
&& =\frac{d}{dh} L(m+1)L(-i_1)\ldots L(-i_k)v. \nonumber
\eea
\end{proof}

After some calculation we can show that (as we mentioned before)
the right  radical of the Shapovalov form for $M_2(c,h)$ coincide with
the radical of Shapovalov form for $M(c,h)$ on the graded subspaces
of degrees one and two.
But on degree three subspace the situation changes.
\begin{example}
Because of Lemma \ref{stupid} we have:
\bea
&& [S_2(h,c)_3]= \left( \begin{array}{cccccc} 
24h(h + 1)(1 + 2h) &  36h(h + 1), & 24h \\
36h(h + 1) & (h + 2)(8h + c) + 18h & 16h+ 2c \\  
 24h & 16h + 2c &  6h + 2c &  \\
0 & 0 & 0 \\
0 & 0 & 0 \\
0 & 0 & 0 
\end{array} \right. \nn
&& \left. \begin{array}{cccccc}  
144h^2+144h+24 & 72h + 36 & 24  \\   
72h + 36 & 16h + 34 + c & 16 \\
24 & 16  & 6 \\ 
24h(h + 1)(1 + 2\,h)  & 36h(h + 1) & 24h \\
36h(h + 1) & (h + 2)(8h + c) + 18h & 16h + 2c \\ 
24h & 16h + 2c & 6h + 2c 
 \end{array}\right).
\eea
Then $$ S_2(h,c)_3=48^2 h^{4}(16h^{2} + 2hc - 10h + c)^2(3h^{2} + hc - 
7h + 2 + c)^2.$$
By analyzing the null space of $[S_2(h,c)_3]$ we see 
that for some pairs $(c,h)$ there are singular vectors in $M_2(c,h)(3)$
that are not contained in $U(Vir)_-.v$. 

For instance, If we let $c=h=1$ then $M'_2(1,1)(3)$ is two dimensional.
Moreover this space is spanned by singular vectors
$$ L(-1)^3v-4L(-1)L(-2)v+6L(-3)v$$
and
$$ -2L(-1)L(-2)v+5L(-3)v+L(-1)^3w-4L(-1)L(-2)w
+6L(-3)w .$$
Therefore 
$${\rm Hom}_{Vir}(M_2(1,4),M_2(1,1)) \cong \mathbb{C}.$$ 
\end{example}
\begin{remark}
After we finished the appendix, we realized that in
\cite{F} and \cite{MRS} a similar problem has been studied. In
particular, in \cite{F} many explicit expressions for
singular vectors inside $M_n(c,h)$ have been obtained.
\end{remark}

\end{document}